\documentclass[11pt,a4paper,fleqn]{amsart}
\usepackage{amssymb}
\usepackage{hyperref}
\usepackage{natbib}

\newcommand{\ee}{\mathbb{E}}

\newcommand{\pp}{\mathbb{P}}
\newcommand{\rr}{\mathbb{R}}
\DeclareMathOperator{\Var}{Var}
\DeclareMathOperator{\sgn}{sgn}
\newcommand{\cn}{\mathcal{N}}
\newcommand{\cl}{\mathcal{L}}
\newcommand{\gs}{S_\alpha(\sigma,\beta,\mu)}

\newcommand{\dto}{\stackrel{d}{\to}}
\newcommand{\de}{\stackrel{d}{=}}
\newcommand{\toi}{\to\infty}
\newcommand{\as}[1]{\quad\text{as }#1\toi}
\newtheorem{thm}{Theorem}
\newtheorem{lem}[thm]{Lemma}
\theoremstyle{remark}
\newtheorem{rmk}{Remark }

\newcommand{\beq}{\begin{equation}}
\newcommand{\eeq}{\end{equation}}
\newcommand{\bt}{\begin{thm}}
\newcommand{\et}{\end{thm}}
\newcommand{\bl}{\begin{lem}}
\newcommand{\el}{\end{lem}}

\begin{document}
\bibliographystyle{plainnat}
\setcitestyle{numbers}
\title{On the functional limits for sums of a function of partial sums}
\author{Kamil Marcin Kosi\'nski}
\thanks{This research was done while the author stayed at the Vrije Universiteit, Amsterdam.}
\email{k.m.kosinski@uva.nl}
\address{Wydzia{\l} Matematyki, Informatyki i Mechaniki, Uniwersytet Warszawski, Warszawa, Poland}
\curraddr{Korteweg-de Vries Institute for Mathematics, University of Amsterdam, P.O. Box 94248,
1090 GE Amsterdam, The Netherlands}

\date{March 6, 2009}
\subjclass[2000]{Primary G0F05}

\begin{abstract}
We derive a functional central limit theorem (fclt) for normalised sums of a function of the partial sums of independent and
identically distributed random variables. In particular, we show, using a technique presented in 
 Huang and Zhang (Electron. Comm. Probab. 12 (2007), 51--56), that the result from 
 Qi (Statist. Probab. Lett. 62 (2003), 93--100), for normalised products of partial sums, can be generalised in this fashion to a fclt.
\end{abstract}

\maketitle

\section{Introduction}
While considering limiting properties of sums of records, \citet{Arnold} obtained the following
version of the central limit theorem (clt) for a sequence $(X_n)$ of independent, identically
distributed (iid) exponential random variables (rv's) with the mean equal one:
\[
  \left(\prod_{k=1}^n\frac{S_k}{k}\right)^{1/\sqrt n}\dto e^{\sqrt 2\cn}\as n,
\]
where $S_n=\sum_{k=1}^n X_k$ and $\cn$ is a standard normal random variable.
\par Later \citet{Remiweso} extended such a clt to general iid positive rv's $(X_n)$. Namely,
provided that $\ee X_1^2<\infty$,
\beq
\label{res1}
  \left(\prod_{k=1}^n\frac{S_k}{k\mu}\right)^{\gamma/\sqrt n}\dto e^{\sqrt 2\cn}\as n,
\eeq
where $\mu=\ee X_1$ and $\gamma=\mu/\sigma$ with $\sigma^2=\Var X_1>0$.
\par This result was generalised by \citet{Qi} by assuming that the underlying distribution of $X_1$
is in the domain of attraction of a stable law with index $\alpha\in(1,2]$. 
In this case
\beq
\label{res2}
  \left(\prod_{k=1}^n\frac{S_k}{k\mu}\right)^{\mu/a_n}\dto e^{(\Gamma(\alpha+1))^{1/\alpha}\cl}\as n,
\eeq
where $\Gamma(\alpha+1)=\int_0^\infty x^\alpha e^{-x}\, dx$ and the sequence $a_n$ is taken such that
\[
  \frac{S_n-n\mu}{a_n}\dto\cl,
\]
where $\cl$ is one of the stable distributions with index $\alpha\in(1,2]$. \citet{LuiQi} obtained a similar 
result in the case $\alpha=1$ with $\ee |X_1|<\infty$. In a paper by \citet{Huang} it is shown that \eqref{res1} 
follows from the weak invariance principle and the whole result can be reformulated to a functional theorem. 
\par Let $(S_n)$ be \textit{any} nondecreasing sequence of positive rv's (which do not have
to be a sequence of partial sums). Suppose there exists a standard
Wiener process $(W(t))_{t\ge0}$ and two positive constants $\mu$ and $\sigma$ such that
\beq
\label{cond1}
  \frac{S_{[nt]}-[nt]\mu}{\sigma\sqrt{n}}\dto W(t)\quad\text{in}\quad D[0,1]\as n,
\eeq
and
\beq
\label{cond2}
  \sup_n \frac{\ee |S_n-n\mu|}{\sqrt n}<\infty.
\eeq
Then the result from \citet{Huang} states that
\beq
\label{res3}
  \left(\prod_{k=1}^{[nt]}\frac{S_k}{k\mu}\right)^{\gamma/\sqrt n}\dto\exp\left\{\int_0^t\frac{W(x)}{x}\,dx\right\}
\quad\text{in}\quad D[0,1]\as n,
\eeq
where $\gamma=\mu/\sigma$. 
\par For example, if $(X_n)$ are iid positive rv's with mean $\mu$ and variance
$\sigma^2$ and $S_n$ are the partial sums, then \eqref{cond1} is satisfied by the invariance principle, cf.
\cite[Theorem 14.1]{Billingsley}, and \eqref{cond2} follows from the Cauchy-Schwarz inequality. Then, one can
check that \eqref{res3} implies \eqref{res1}.
\par The purpose of this paper is to show that the technique used by \citet{Huang} can be utilised to obtain a similar
result for the rv's in the domain of attraction of a stable law with index $\alpha\in(1,2]$. Furthermore, we will
set our discussion in a more general setting. It is straightforward that \eqref{res1}, and analogously \eqref{res2}, is
a simple corollary from
\[
  \frac{\sum_{k=1}^n f(S_k/k) - b_n}{a_n}\dto\cn,
\]
if one sets $f(x)=\log x$ and chooses the sequences $a_n$, $b_n$ properly.
\section{Preliminaries}
The following theorem is well known and can be easily found in the literature, see e.g., \cite{Bingham}.
\bt[Stability Theorem]
\label{general:stab}
  The general stable law with index $\alpha\in(0,2]$ is given by a characteristic function of one of the following forms:
\begin{enumerate}
  \item $\phi(t)=\exp(-\sigma^\alpha|t|^\alpha(1-i\beta(\sgn t)\tan\frac{1}{2}\pi\alpha)+i\mu t)$, $\alpha\ne1$,
  \item $\phi(t)=\exp(-\sigma|t|(1+i\beta(\sgn t)\frac{2}{\pi}\log|t|)+i\mu t)$,  $\alpha=1$,
  \item $\phi(t)=\exp(-\sigma^2t^2/2+i\mu t)$, $\alpha=2$
\end{enumerate}
with $\beta\in[-1,1]$, $\mu\in\rr$ and $\sigma>0$.
\et
From the above theorem one can see that every stable law with index $\alpha\in(0,2)$ can be parametrized by four parameters and written
as $\gs$. We distinguish the case $\alpha=2$ because otherwise $S_2(\sigma,\beta,\mu)\de\cn(\mu,2\sigma^2)$ and $\beta$
plays no role, moreover one would like to think of $S_2(1,\beta,0)$ as $\cn$ not $\cn(0,2)$. 
\par Let $(X_n)$ be a sequence of iid rv's, set $S_n=\sum_{k=1}^n X_k$ and assume $X_1$ is in the
domain of attraction of a stable law with index $\alpha\in(1,2]$. Note that for such $X_1$ we have $\ee |X_1|<\infty$.
Recall that a sequence of iid rv's $(X_n)$ is said to be in the domain of attraction of a stable law $\gs$, if
there exists constants $a_n>0$ and $b_n\in\rr$ such that
\beq
\label{def}
  \frac{S_n-b_n}{a_n}\dto\gs.
\eeq
Clearly by scaling it suffices to let $\sigma=1$ and $\mu=0$, hence only the parameters $\alpha$ and $\beta$ are
unaltered by scaling. In case $\alpha=2$ as mentioned above, we understand $S_2(1,\beta,0)$ as $\cn$.
Moreover, in case $\alpha\in(1,2]$, the sequence $b_n$ can be taken equal to $\mu n$.
On how to choose $a_n$ we refer to \cite[Theorem 4.5.1]{Whitt}. The choice is irrelevant for our discussion,
the only fact that plays the crucial role is the fact that if \eqref{def} holds, then $a_n=n^{1/\alpha}L(n)$, where
$L$ is slowly varying.
\par Furthermore, in addition to stable clt \eqref{def}, there is convergence in distribution
\beq
\label{stable}
  S_n(t):=\frac{S_{[nt]}-[nt]\mu }{a_n}\dto \cl(t)\quad\text{in}\quad D[0,1],
\eeq
where $\cl$ is a standard $(\alpha,\beta)$-stable L\'evy motion, with
\[
  \cl(t)\de t^{1/\alpha}S_\alpha(1,\beta,0)\de S_\alpha(t,\beta,0)
\]
and as before, for $\alpha=2$, $\cl$ is a standard Wiener process, cf. \cite[Theorem 4.5.3]{Whitt}.
\section{Main Results}
\bt
  \label{result}
Let $(X_n)$ be a sequence of iid rv's in the domain of attraction of
the stable law $S_\alpha(1,\beta,0)$ with $\alpha\in(1,2]$, so that \eqref{def} holds
for some sequence $a_n$ and $b_n=n\mu$, where $\mu=\ee X_1$. Let $f$ be a real function defined on an interval $I$
such that $\pp(X_1\in I)=1$ and $f'(\mu)$ exists. Then, as $n\toi$
\beq
  \label{tw3:teza}
\frac{1}{a_n}\sum_{k=1}^{[nt]}\left(f(S_k/k)-f(\mu)\right)\dto f'(\mu)\int_0^t\frac{\cl(x)}{x}\,dx\quad\text{in}\quad D[0,1],
\eeq
where $S_k$ denotes the $k$-th partial sum.
\et
Because $\cl(x)$ is c\'adl\'ag, it has at most countably many discontinuity points, so the integral 
on the right hand side of \eqref{tw3:teza} exists and is finite almost surely if
\[
  \int_0^1\frac{|\cl(x)|}{x}\,dx<\infty\quad\text{a.s}.
\]
To ensure this, note that for a positive nondecreasing function $h$ we have
\[
  \int_0^t\frac{|\cl(x)|}{x}\,dx\le\sup_{0\le s\le t}\left\{\frac{|\cl(s)|}{h(s)}\right\}\int_0^t\frac{h(x)}{x}\, dx.
\]
Setting $h(x)=x^\gamma$ with $\gamma\in(0,1/\alpha)$ we get $\int_0^1\frac{h(x)}{x}\, dx<\infty$ and
\[
  \sup_{0\le s\le t}\left\{\frac{|\cl(s)|}{h(s)}\right\}\to 0\,\text{ a.s. as}\quad t\to 0,
\]
by Khintchine's Theorem, see e.g., \cite[Theorem 2.1]{Nielsen}. This guaranties the existence of the integral in 
\eqref{tw3:teza} as well as implies that
\[
  \sup_{0\le s\le t}\left|\int_0^t\frac{\cl(x)}{x}\,dx\right|=0\,\text{ a.s. as}\quad t\to 0,
\]
a fact that is going to be used later in the proof.
\begin{rmk}
\label{remark:calc}
Observe that
\begin{align*}
  \int_0^t\frac{\cl(x)}{x}\,dx &=\lim_{n\toi}\sum_{k=1}^n\frac{t}{n}\frac{n}{tk}\cl\left(\frac{tk}{n}\right)
=\lim_{n\toi}\sum_{k=1}^n\sum_{i=1}^k\frac{1}{k}\left(\cl\left(\frac{ti}{n}\right)-\cl\left(\frac{t(i-1)}{n}\right)\right)\\
&=\lim_{n\toi}\sum_{i=1}^n\left(\cl\left(\frac{ti}{n}\right)-\cl\left(\frac{t(i-1)}{n}\right)\right)\sum_{k=i}^n\frac{1}{k}\\
&\de \lim_{n\toi}\left(\sum_{i=1}^n\left(\sum_{k=i}^n\frac{1}{k}\right)^\alpha\right)^{1/\alpha}S_\alpha\left(\frac{t}{n},\beta,0\right)\\
&\de  S_\alpha\left(t,\beta,0\right) \lim_{n\toi}\left(\sum_{i=1}^n\frac{1}{n}\left(\sum_{k=i}^n\frac{1}{k}\right)^\alpha\right)^{1/\alpha}\\
&=S_\alpha\left(t,\beta,0\right)\left(\int_0^1(-\log x)^\alpha\,dx\right)^{1/\alpha}\\
&=S_\alpha\left(t,\beta,0\right)\left(\Gamma(\alpha+1)\right)^{1/\alpha}.
\end{align*}
If $X_1$ is a positive rv, then the limiting stable law has $\beta=1$. Setting 
$f(x)=\mu\log(x/\mu)$, Theorem \ref{result} yields
\[
  \left(\prod_{k=1}^n\frac{S_k}{k\mu}\right)^{\mu/a_n}\dto
 \exp\left(S_\alpha\left(\Gamma(\alpha+1),1,0\right)\right),\as n.
\]
which is the result \eqref{res2} obtained by \citet{Qi}.
\end{rmk}
\begin{rmk}
If $\ee X_1^2<\infty$, then $X_1$ is in the domain of attraction of normal distribution $\cn$
and $a_n\sim\sigma\sqrt n$, where $\sigma^2=\Var(X_1)>0$. If furthermore $X_1$ is positive, then setting $\gamma=\mu/\sigma$
\[
  \left(\prod_{k=1}^{[nt]}\frac{S_k}{k\mu}\right)^{\gamma/\sqrt n}\dto \exp\left(\int_0^t\frac{W(x)}{x}\,dx\right),\quad\text{in}\quad D[0,1]\as n,
\]
which coincides with the result \eqref{res3} by \citet{Huang}.
\end{rmk}
Before proceeding to the proof of the main theorem, we need a technical lemma.
\bl Under the assumptions of Theorem \ref{result} 
\label{lemma:sb}
\[
  \sum_{k=1}^n\frac{\ee|S_k-k\mu|}{k}=O(a_n).
\]
\el
\begin{proof}
Note that
\[
  \sum_{k=1}^n\frac{\ee|S_k-k\mu|}{k} \le
\sup_{k\le n}\left\{\ee\frac{|S_k-k\mu|}{a_k}\right\}\sum_{k=1}^n\frac{a_k}{k}.
\]
By Theorem 6.1 in \citet{De}
\beq
\label{cond:sim2}
  \ee\frac{|S_n-n\mu|}{a_n}=O(1).
\eeq
Now, for a regularly varying function $A>0$ with index $\gamma>-1$, its easy to see that
\[
  \sum_{k\le x} A(k)\sim\int_1^x A(t)dt\sim\frac{1}{1+\gamma}xA(x)\as x\quad\text{if }\gamma>-1,
\]
where the last asymptotic equivalence follows from the Karamata's Theorem, cf. \cite[Theorem 1.5.8]{Bingham}.
Recall that $a_n$ is slowly varying with index $1/\alpha>0$, this implies  
\[
  \sum_{k=1}^n\frac{a_k}{k}=O(a_n),
\]
and proves the Lemma.
\end{proof}
Now we may proceed to the proof of the main theorem. The proof follows the steps of the proof of \eqref{res3} in
\citet{Huang}.
\begin{proof}[Proof of Theorem \ref{result}]
 Expand $f$ in the neighbourhood of $\mu$, then
\beq
\label{expansion}
\frac{1}{a_n}\sum_{k=1}^{[nt]}(f(S_k/k)-f(\mu))=\frac{f'(\mu)}{a_n}\sum_{k=1}^{[nt]}(S_k/k-\mu)+\frac{1}{a_n}\sum_{k=1}^{[nt]}(S_k/k-\mu)r(S_k/k),
\eeq
where $r(x)\to 0$ as $x\to\mu$. Note that $\ee |X_1|<\infty$ so by the SLLN $r(S_k/k)\to 0$ a.s..
It now follows from Lemma \ref{lemma:sb} that
\[
  \sup_{0\le t\le 1}\left|\frac{1}{a_n}\sum_{k=1}^{[nt]}(S_k/k-\mu)r(S_k/k)\right|
\le \frac{1}{a_n}\sum_{k=1}^{n}\frac{|S_k-k\mu|}{k}|r(S_k/k)|=o_\pp(1).
\]
So, according to \eqref{expansion} it suffices to show that, as $n\toi$
\beq
\label{cond}
  Y_n(t):=\frac{1}{a_n}\sum_{k=1}^{[nt]}\frac{S_k-k\mu}{k}\dto\int_0^t\frac{\cl(x)}{x}\,dx,\quad\text{in}\quad D[0,1].
\eeq
Let
\[
  H_\varepsilon(f)(t)=\left\{
 \begin{array}{rl}
 \int_\varepsilon^t\frac{f(x)}{x}\,dx,& \varepsilon<t\le 1\\
 0,& 0\le t\le\varepsilon
 \end{array}\right.
\]
and
\[
  Y_{n,\varepsilon}(t)=\left\{
 \begin{array}{rl}
 \frac{1}{a_n}\sum_{k=[n\varepsilon]+1}^{[nt]}\frac{S_k-k\mu}{k},& \varepsilon<t\le 1\\
 0,& 0\le t\le\varepsilon,
 \end{array}\right.
\]
It is obvious that
\beq
\label{eq:1}
  \sup_{0\le t\le 1}\left|\int_{0}^t\frac{\cl(x)}{x}\,dx-H_\varepsilon(\cl)(t)\right|=
\sup_{0\le t\le \varepsilon}\left|\int_{0}^t\frac{\cl(x)}{x}\,dx\right|\to0\,\text{ a.s. as }\,\varepsilon\to 0
\eeq
and
\beq
\label{eq:2}
  \ee\max_{0\le t\le 1}\left|Y_n(t)-Y_{n,\varepsilon}(t)\right|\le 
\frac{1}{a_n}\sum_{k=1}^{[n\varepsilon]}\frac{\ee|S_k-k\mu|}{k}\le C \varepsilon^{1/\alpha}
\eeq
by the same argumentation as in the proof of Lemma \ref{lemma:sb}. 
\par On the other hand, it is
easily seen that, for $n$ large enough such that $n\varepsilon\ge 1$,
\begin{align*}
\sup_{\varepsilon\le t\le 1}\Bigg|\sum_{k=[n\varepsilon]+1}^{[nt]}\frac{S_k-k\mu}{k}&-
\int_{n\varepsilon}^{nt}\frac{S_{[x]}-[x]\mu}{x}\,dx\Bigg|
=\sup_{\varepsilon\le t\le 1}\left|\int_{[n\varepsilon]+1}^{[nt]+1}\frac{S_{[x]}-[x]\mu}{[x]}\,dx-
\int_{n\varepsilon}^{nt}\frac{S_{[x]}-[x]\mu}{x}\,dx\right|\\
&\le\left |\int_{n\varepsilon}^{[n\varepsilon]+1}\frac{S_{[x]}-[x]\mu}{x}\,dx \right|
+\sup_{\varepsilon\le t\le 1}\left|\int_{nt}^{[nt]+1}\frac{S_{[x]}-[x]\mu}{x}\,dx\right|\\
&\quad+\sup_{\varepsilon\le t\le 1}\left|\int_{[n\varepsilon]+1}^{[nt]+1}(S_{[x]}-[x]\mu)\left(\frac{1}{x}-\frac{1}{[x]}\right)\,dx\right|\\
&\le\max_{k\le n}|S_k-k\mu|\sup_{\varepsilon\le t\le 1}\left(\frac{2}{n\varepsilon}+\frac{2}{nt}+\frac{1}{n\varepsilon}\right)\\
&\le 5\max_{k\le n}|S_k-k\mu|/(n\varepsilon)=O_\pp(a_n/n)=o_\pp(1),
\end{align*}
by noticing that $\max_{k\le n} |S_k-k\mu|/a_n\dto \sup_{0\le t\le 1} |\cl(t)|$ according to \eqref{stable}.
So
\[
  \frac{1}{a_n}\sum_{k=[n\varepsilon]+1}^{[nt]}\frac{S_k-k\mu}{k}=\frac{1}{a_n}\int_{n\varepsilon}^{nt}\frac{S_{[x]}-[x]\mu}{x}\,dx+o_\pp(1)
=\int_{\varepsilon}^{t}\frac{S_{n}(x)}{x}\,dx+o_\pp(1)
\]
uniformly in $t\in[\varepsilon,1]$. Notice that $H_\varepsilon(\cdot)$ is a continuous mapping on the space $D[0,1]$. Using
the continuous mapping theorem (cf. \cite[Theorem 2.7]{Billingsley}) it follows that
\beq
\label{eq:3}
  Y_{n,\varepsilon}(t)=H_\varepsilon(S_n)(t)+o_\pp(1)\dto H_\varepsilon(\cl)(t)\quad\text{in}\quad D[0,1]\as n.
\eeq
Combining \eqref{eq:1}-\eqref{eq:3} yields \eqref{cond} by \cite[Theorem 3.2]{Billingsley}.
\end{proof}
\section{Extensions}
To prove Lemma \ref{lemma:sb}, we have only used the property \eqref{cond:sim2} (which is in fact
the condition \eqref{cond2}) and the fact that $a_n$ varies regularly with a positive index. The proof
of Theorem \ref{result} was based on the convergence \eqref{stable} and the fact that $S_k/k\to\mu$ a.s.. All 
those conditions are satisfied when $S_k$ is defined to be the partial sum of a sequence of iid rv's in the
domain of attraction of a stable law with index greater than one. However, we do not need to assume anything
about $S_k$ and only require that it satisfies the aforementioned conditions. This leads to
\bt
\label{result2}
Let $(S_k)$ be a sequence of random variables. Suppose there exists an $(\alpha,\beta)$-stable L\'evy process $(\cl(t))_{t\ge 0}$,
a constant $\mu$ and a sequence $a_n$ such that as $n\toi$
\[
  \frac{S_{[nt]}-[nt]\mu}{a_n}\dto\cl(t)\quad\text{in}\quad D[0,1],
\]
where $a_n$ can be written as $a_n=n^{1/\alpha}L(n)$ with $\alpha\in(1,2]$ and $L$ slowly varying. In addition,
suppose that
\beq
\label{cond2:again}
  \sup_n \frac{\ee |S_n-n\mu|}{a_n}=O(1),
\eeq
and $S_n/n\to\mu$ a.s., then, as $n\toi$
\[
\frac{1}{a_n}\sum_{k=1}^{[nt]}\left(f(S_k/k)-f(\mu)\right)\dto f'(\mu)\int_0^t\frac{\cl(x)}{x}\,dx\quad\text{in}\quad D[0,1],
\]
for any real function $f$ defined on an interval $I$ such that $\pp(S_k/k\in I)=1$ for all $k$, provided that $f'(\mu)$ exists.
\et
In their paper, Huang and Zhang showed that if $(S_k)$ is a nondecreasing (in fact we only need monotonicity) sequence
satisfying \eqref{cond2:again}, then $S_k/k\to\mu$ a.s.. Thus, Theorem \ref{result2} is an extension of the result \eqref{res3}
from \citet{Huang}.
\section{Acknowledgments}
The author is grateful to dr. R. Lata{\l}a for supervision of his work and would also like to thank the referees for 
pointing out some subtleties in the previous version of this paper.

\end{document}